\documentclass[]{article}
\begin{document}
\newtheorem{proposition}{Proposition}[section]
\newtheorem{definition}{Definition}[section]
\newtheorem{lemma}{Lemma}[section]

\title{\bf  Calculus on Dual Real Numbers}
\author{Keqin Liu\\Department of Mathematics\\The University of British Columbia\\Vancouver, BC\\
Canada, V6T 1Z2}
\date{Augest, 2018}
\maketitle

\begin{abstract} We present the basic theory of calculus on dual real numbers, and prove the counterpart of the odinary fundamental theorem of calculus in the context of dual real numbers.
\end{abstract}

\medskip
The purpose of this paper is to study calculus on dual real numbers. Unlike the multi-variables calculus on the Cartesian product of finite many copies of the real number field and the complex analysis on the complex number field, the generalizations of the order relation on the real number field plays a central role in the theory of calculus on dual real numbers. Hence, calculus on dual real numbers seemd to be closer to the well-known single variable calculus than both multi-variables calculus and complex analysis.

\medskip
The main result of this paper is to explain how to develpop the basic theory of calculus on dual real numbers. In section 1, we make the dual real number algebra into a normed algebra and introduce two generalizations of the order relatio on the real number field. In section 2, we define the differentiability in dual real numbers, and characterize the differentiability by using the real-valued component functions of a dual real number-valued function. In section 3, we introduce two types of integrals based on the two generalized order relations and prove the counterpart of the odinary fundamental theorem of calculus in the context of dual real numbers.

\medskip
\section{Two generalized order relations}

\medskip
We begin this section by recalling some facts about dual real numbers. For any two elements $(x_1,\, x_2)$  and $(y_1,\, y_2)$ from the $2$-dimensional real vector space 
$\mathcal{R}^2=\mathcal{R}\times\mathcal{R}$, we  define their product according to the the fololowing rule called 
{\bf dual number multiplication}:
$$(x_1,\, x_2)\,\cdot\,(y_1, \, y_2):= (x_1y_1,\, x_1y_2+x_2y_1).$$
The vector space $\mathcal{R}^2$ with respect to the dual number multiplication is a real associative algebra, which is called the  {\bf dual real number algebra} and denoted by $\mathcal{R}^{(2)}$. An element of $\mathcal{R}^2$ is called a dual real number. Clearly, the  dual real number algebra 
$\mathcal{R}^{(2)}$ is both unital and commutative. We denote the multiplication identity $(1,\, 0)$ by $1$, and the element
$(0,\, 1)$ by $1^{\#}$. Then every dual real number  $(x_1,\, x_2)$ of $\mathcal{R}^{(2)}$ can be expressed in a unique way as a linear combination of $1$ and $1^{\#}$:
$$ x=(x_1,\, x_2)=x_1+x_2\,1^{\#}\qquad \mbox {for $x_1$, $x_2\in \mathcal{R}^{(2)}$},$$
where $Re\,x:=x_1$ and $Ze\,x:=x_2$ are called the {\bf real part} and the {\bf zero-divisor part} of $x$, respectively.

\medskip
The dual real number algebra $\mathcal{R}^{(2)}$ is not a field and has many zero-divisors.
In fact, if $0\ne x\in \mathcal{R}^2$, then $x$ is a zero-divisor if and only if $Re\,x=0$ and
$x$ is invertible if and only if $Re\,x\ne 0$. Moreover if $x$ is invertible,  then  the inverse $x^{-1}$ of $x$ is given by
$x^{-1}=\displaystyle\frac{1}{Re\, x}-\displaystyle\frac{Ze\, x}{(Re\, x)^2}\,1^{\#} .$

\medskip
\begin{definition}\label{def1.1} The real-valued function $||\quad||: \mathcal{R}^{(2)}\to \mathcal{R}$  defined by
\begin{equation}\label{eq1}
||x||:=\sqrt{2\,(Re\,x)^2+(Ze\,x)^2}\quad\mbox{for $x\in\mathcal{R}^{(2)}$}
\end{equation}
is called the {\bf norm} in $\mathcal{R}^{(2)}$
\end{definition}

\medskip
The dual real number algebra is a normed algebra with respect to the norm introduced in Definition \ref{def1.1} by the following proposition.

\medskip
\begin{proposition}\label{pr1.1} Let $x$, $y\in\mathcal{R}^{(2)}$ and $a\in\mathcal{R}$.
\begin{description} 
\item[(i)] $||x||\ge 0$, with equality only when $x=0$.
\item[(ii)] $||ax||=|a|\,||x||$, where $|a|$ denotes the absolute value of the real number $a$.
\item[(iii)]  $||x+y||\le  ||x||+||y||$.
\item[(iv)]  $||xy||\le  ||x||\,||y||$.
\end{description}
\end{proposition}

\medskip
\noindent
{\bf Proof} a direct computation.

\hfill\raisebox{1mm}{\framebox[2mm]{}}

\medskip
Unlike the Cartesian product of finite many copies of the real number field and the complex field, there are two generalized order relations on $\mathcal{R}^{(2)}$  which are compatible with the multiplication in $\mathcal{R}^{(2)}$.

\begin{definition}\label{def1.2} Let $x$ and $y$ be two elements of $\mathcal{R}^{(2)}$.
\begin{description}
\item[(i)] We say that $x$ is {\bf type 1 greater than} $y$ ( or $y$ is {\bf type 1 less than} $x$) and we write $x\stackrel{1}{>} y$ (or $y\stackrel{1}{<} x$) if
$$
\mbox{either}\quad \left\{\begin{array}{c}
Re\,x>Re\,y\\ Ze\,x\ge Ze\,y\end{array}\right.\quad\mbox{or}\quad
\left\{\begin{array}{c}
Re\,x=Re\,y\\ Ze\,x> Ze\,y\end{array}\right.
$$
\item[(ii)] We say that $x$ is {\bf type 2 greater than} $y$ ( or $y$ is {\bf type 2 less than} $x$) and we write $x\stackrel{2}{>} y$ (or $y\stackrel{2}{<} x$) if
$$
\mbox{either}\quad \left\{\begin{array}{c}
Re\,x>Re\,y\\ Ze\,y\ge Ze\,x\end{array}\right.\quad\mbox{or}\quad
\left\{\begin{array}{c}
Re\,x=Re\,y\\ Ze\,y> Ze\,x\end{array}\right.
$$
\end{description}
\end{definition}

We use $x\stackrel{\theta}{\ge} y$ when $x\stackrel{\theta}{>} y$ or $x=y$ for $\theta =1,\, 2$.
By Definition 1.1, if $Re\,x=Re\,y$, then $x\stackrel{1}{>} y\Longleftrightarrow y\stackrel{2}{>} x$; if $Ze\,x=Ze\,y$, then $x\stackrel{1}{>} y\Longleftrightarrow x\stackrel{2}{>} y$.
The following proposition gives the basic properties of the two generalized order relations.

\begin{proposition}\label{pr1.2} Let $x$, $y$ and $z$ be elements of $\mathcal{R}^{(2)}$ and $\theta =1,\, 2$.
\begin{description}
\item[(i)] One of the following holds:
$$ x\stackrel{1}{>} y, \quad y\stackrel{1}{>} x, \quad x=y, \quad x\stackrel{2}{>} y, \quad y\stackrel{2}{>} x .$$
\item[(ii)] If $x\stackrel{\theta}{>} y$ and $y\stackrel{\theta}{>} z$, then $x\stackrel{\theta}{>} z$.
\item[(iii)] If $x\stackrel{\theta}{>} y$, then $x+z\stackrel{\theta}{>} y+z$.
\item[(iv)] If $x\stackrel{\theta}{>} 0$ and $y\stackrel{\theta}{>} 0$, then $xy\stackrel{\theta}{\ge} 0$.
\item[(v)] If $x\stackrel{\theta}{>} y$, then $-x\stackrel{\theta}{<} -y$.
\end{description}
\end{proposition}

\medskip
\noindent
{\bf Proof} Clear.

\hfill\raisebox{1mm}{\framebox[2mm]{}}

\medskip
\section{Differentiation}

\medskip
By Proposition \ref{pr1.1},  the dual real number algebra $\mathcal{R}^{(2)}$ is a metric space with the distance function $||\,\,||$.  If $c\in\mathcal{R}^{(2)}$ and $\epsilon\in \mathcal{R}$, we use $N(c;\,\epsilon)$  and $N^*(c;\,\epsilon)$ to denote the ordinary $\epsilon$-{\bf neighborhood} and {\bf deleted $\epsilon$-neighborhood} of $c$, respectively, i.e.,  
$$N(c;\,\epsilon)=\{x\in \mathcal{R}^{(2)}\,|\, ||x-c||<\epsilon\}\quad\mbox{and}\quad
N^*(c;\,\epsilon):=N(c;\,\epsilon)\setminus\{c\}.$$
For $\theta\in\{1,\,2\}$, the  set
$
N_{\theta}(c;\,\epsilon):=\{x\in N(c;\,\epsilon)\,\big|\, \mbox{$x\stackrel{\theta}{\ge}c$ or
 $x\stackrel{\theta}{\le}c$}\}
$
is called the {\bf type $\theta$ $\epsilon$- neighborhood} of $c$ and the set
$
N_{\theta}^*(c;\,\epsilon):=N_{\theta}(c;\,\epsilon)\setminus\{c\}
$
is called the {\bf deleted type $\theta$ $\epsilon$- neighborhood} of $c$. 
By Proposition\ref{pr1.2} (i), we have
$$
N(c;\,\epsilon)=N_{1}^*(c;\,\epsilon)\bigcup N_{2}^*(c;\,\epsilon).
$$

\medskip
We now introduce the differentiability in the following

\medskip
\begin{definition}\label{def2.1} Let $D$ be an open subset of $\mathcal{R}^{(2)}$ and let 
$c\in D$. 
\begin{description} 
\item[(i)] We say that $f: D\to \mathcal{R}^{(2)}$ is {\bf type $\theta$ differentiable} at $c$ with $\theta\in\{1,\, 2\}$ if for each positive real number $\epsilon>0$ there exist a  positive real number $\delta>0$ and a dual real number $f_{\theta}'(c)\in \mathcal{R}^{(2)}$ such that
$$ x\in N_{\theta}^*(c;\,\delta)\subseteq D\Rightarrow 
\frac{f(x)-f(c)-f_{\theta}'(c)(x-c)}{||x-c||}\in N(0;\,\epsilon).  $$
The dual real number $f_{\theta}'(c)$ is called the {\bf type $\theta$ derivative} of $f$ at $c$, which is also denoted by $\displaystyle\frac{df}{d_{\theta}x}(c)$. 
\item[(ii)] We say that $f: D\to \mathcal{R}^{(2)}$ is {\bf differentiable} at $c$ if for each positive real number $\epsilon>0$ there exist a  positive real number $\delta>0$ and a dual real number $f'(c)\in \mathcal{R}^{(2)}$ such that
$$ x\in N^*(c;\,\delta)\subseteq D\Rightarrow 
\frac{f(x)-f(c)-f'(c)(x-c)}{||x-c||}\in N(0;\,\epsilon).  $$
The dual real number $f'(c)$ is called the {\bf derivative} of $f$ at $c$, which is also denoted by $\displaystyle\frac{df}{dx}(c)$. If $f$ is  differentiable at each point of the open subset $D$, then $f$ is said to be {\bf differentiable} on $D$.
\end{description}
\end{definition}

\medskip
It is easy to check that if a dual real number-valued $f$ is  type $\theta$ differentiable at 
$c\in \mathcal{R}^{(2)}$, then the type $\theta$ derivative of $f$ at $c$ is unique for $\theta\in\{1,\, 2\}$.

\bigskip
Let $S$ be a subset of $\mathcal{R}^{(2)}$. A function $f: S\to \mathcal{R}^{(2)}$ can be expressed as  
$$f(x)=u(x)+v(x)\,1^{\#}\quad\mbox{for $x=x_1+x_2\,1^{\#}\in \mathcal{R}^{(2)}$},$$
where  $u(x):=u(x_1,\,x_2)$ and $v(x):=v(x_1,\,x_2)$ are two real-valued functions of two real variables $x_1$ and $x_2$, which are called the {\bf real component} and the {\bf zero-divisor component} of $f$, respectively. The following proposition provides an useful characterization of differentiability for  dual 
real-valued functions in terms of their real and zero-divisor components.

\medskip
\begin{proposition}\label{pr2.1} Let $f: D\to \mathcal{R}^{(2)}$ be a  dual 
real-valued function given by
$$f(x)=u(x_1,\,x_2)+v(x_1,\,x_2)\,1^{\#}\quad\mbox{for $x=x_1+x_2\,1^{\#}\in \mathcal{R}^{(2)}$},$$
where $D$ is an open subset of $\mathcal{R}^{(2)}$,  $u(x_1,\,x_2)$ and  $v(x_1,\,x_2)$ are  the real component and the zero-divisor component of $f$, respectively. Let $c=c_1+c_2\,1^{\#}\in D$ with 
$c_1$, $c_2\in \mathcal{R}$.
\begin{description} 
\item[(i)] If the first-order partial derivatives $u_{x_1}$, $u_{x_2}$, $v_{x_1}$ and $v_{x_2}$ exist at
$(c_1,\,c_2)$ and are continuous at $(c_1,\,c_2)$, and the following equations
\begin{equation}\label{eq2}
u_{x_1}=v_{x_2}, \qquad u_{x_2}=0
\end{equation}
hold at $(c_1,\,c_2)$, then $f$ is differentiable at $c$ and the  derivative $f'(c)$ of $f$ at $c$ is given by
\begin{equation}\label{eq3}
f'(c)=u_{x_1}(c)+v_{x_1}(c)\,1^{\#}.
\end{equation}
\item[(ii)] If $f$ is differentiable at $c$, then the equations in  (\ref{eq2}) hold at  $c=c_1+c_2\,1^{\#}$. In this case,  the  derivative $f'(c)$ of $f$ at $c$ is given by (\ref{eq3}).
\end{description}
\end{proposition}

\medskip
\noindent
{\bf Proof} The proof of Proposition \ref{pr2.1} is similar to the proof of the famous fact which characterizes the complex differentiability by using Cauch-Riemmann equations.

\medskip
Foe example, let us prove (ii). If $f$ is differentiable at $c$, then exists $L=L_1+L_2\,1^{\#}$ with $L_1$, $L_2\in  \mathcal{R}$ such that for every $\epsilon>0$, there exists a $\delta>0$ such that $N(c;\,\delta)\subseteq D$ and
\begin{equation}\label{eq4}
x\in N^*(c;\,\delta)\Rightarrow  
\frac{\|f(x)-f(c)- L\,(x-c)\|}{\|x-c\|}<\epsilon .
\end{equation}

\medskip
By dual number multiplication, we have
\begin{eqnarray}\label{eq5}
&&f(x)-f(c)- L\,(x-c)=u(x_1,\,x_2)-u(c_1,\,c_2)-L_1(x_1-c_1)+\nonumber\\
&&\quad +[v(x_1,\,x_2)-v(c_1,\,c_2)-L_1(x_2-c_2)-L_2(x_1-c_1)]\,1^{\#}.
\end{eqnarray}

\medskip
Let $x_2=c_2$ and choose $x_1$ such that $0<|x_1-c_1|<\displaystyle\frac{\delta}{\sqrt{2}}$. Then
$$0<\|x-c\|=\|x_1-c_1\|=\sqrt{2} |x_1-c_1|<\delta ,$$
which implies that $x=x_1+c_2\,1^{\#}\in  N^*(c;\,\delta)$. By  (\ref{eq4}) and  (\ref{eq5}), we get
\begin{eqnarray}\label{eq6}
\epsilon&>&\frac{\|f(c_1+x_2\,1^{\#})-f(c)- L\,(c_1+x_2\,1^{\#}-c)\|}{\|c_1+x_2\,1^{\#}-c\|}\nonumber\\
&=&\frac{1}{\sqrt{2}}\Big\{2\Big[\frac{u(x_1,\,c_2)-u(c_1,\,c_2)}{x_1-c_1}-L_1\Big]^2+\nonumber\\
&&\quad +\Big[\frac{v(x_1,\,c_2)-v(c_1,\,c_2)}{x_1-c_1}-L_2\Big]^2\Big\}^{\frac12}.
\end{eqnarray}
It follows from  (\ref{eq6}) that for every $\epsilon>0$, there exists a $\delta>0$ such that 
$$
\left|\frac{u(x_1,\,c_2)-u(c_1,\,c_2)}{x_1-c_1}-L_1\right|<\epsilon\quad\mbox{and}\quad
\left|\frac{v(x_1,\,c_2)-v(c_1,\,c_2)}{x_1-c_1}-L_2\right|<\epsilon\sqrt{2}
$$
whenever $0<|x_1-c_1|<\displaystyle\frac{\delta}{\sqrt{2}}$. This proves that
\begin{equation}\label{eq7}
u_{x_1}(c_1,\,c_2)=L_1\quad\mbox{and}\quad v_{x_1}(c_1,\,c_2)=L_2.
\end{equation}

\medskip
Similarly, let $x_1=c_1$ and choose $x_2$ such that $0<|x_2-c_2|<\delta$. Then
$$0<\|x-c\|=\|(x_2-c_2)\,1^{\#}\|= |x_2-c_2|<\delta ,$$
which implies that $x=c_1+x_2\,1^{\#}\in  N^*(c;\,\delta)$. By  (\ref{eq4}) and  (\ref{eq5}), we get
\begin{eqnarray*}
\epsilon&>&\frac{\|f(c_1+x_2\,1^{\#})-f(c)- L\,(c_1+x_2\,1^{\#}-c)\|}{\|c_1+x_2\,1^{\#}-c\|}\nonumber\\
&=&\left\{2\left[\frac{u(c_1,\,x_2)-u(c_1,\,c_2)}{x_2-c_2}\right]^2
+\left[\frac{v(c_1,\,x_2)-v(c_1,\,c_2)}{x_2-c_2}-L_1\right]^2\right\}^{\frac12}.
\end{eqnarray*}
which implies that for every $\epsilon>0$, there exists a $\delta>0$ such that 
$$
\left|\frac{u(c_1,\,x_2)-u(c_1,\,c_2)}{x_2-c_2}\right|<
\frac{\epsilon}{\sqrt{2}
}\quad\mbox{and}\quad
\left|\frac{v(c_1,\,x_2)-v(c_1,\,c_2)}{x_2-c_2}-L_1\right|<\epsilon
$$
whenever $0<|x_2-c_2|<\delta$. This proves that
\begin{equation}\label{eq8}
u_{x_2}(c_1,\,c_2)=0\quad\mbox{and}\quad v_{x_2}(c_1,\,c_2)=L_1.
\end{equation}

By (\ref{eq7}) and (\ref{eq8}), (ii) holds.

\hfill\raisebox{1mm}{\framebox[2mm]{}}

\medskip
\section{Type $\theta$ Integrals}

\medskip
In the remaining of this paper, $\theta$ always denote an element in the set $\{1, 2\}$. 
Let  $f: S\to \mathcal{R}^{(2)}$ be a function on a subset $S$ of  $\mathcal{R}^{(2)}$. For convenience, we will use $f_{Re}$ and $f_{Ze}$ to denote the  real component and the  zero-divisor component of a function $f: S\to \mathcal{R}^{(2)}$, respectively. Thus, we have
$$f(x)=f_{Re}(x)+f_{Ze}(x)\,1^{\#}\quad\mbox{for $x=x_1+x_2\,1^{\#}\in \mathcal{R}^{(2)}$}.$$
We say that the function $f: S\to \mathcal{R}^{(2)}$ is {\bf bounded on $S$} if both $f_{Re}$ and $f_{Ze}$ are bounded on $S$ ($\subseteq \mathcal{R}^2= \mathcal{R}\times  \mathcal{R}$).

\medskip
Let $a$, $b\in \mathcal{R}^{(2)}$ and $a\stackrel{\theta}{<}b$. The {\bf type $\theta$ closed interval} $[a, b]_{\theta}$ is defined by
$$
[a, b]_{\theta}:=\big\{\, x\in \mathcal{R}^{(2)}\,|\, a\stackrel{\theta}{\le}x\stackrel{\theta}{\le}b\,\big\}.
$$

\medskip
\begin{definition}\label{def3.1} Let $[a, b]_{\theta}$ be a type $\theta$ closed integral. A {\bf partition} $P$ of  $[a, b]_{\theta}$ is a finite set of points $\{x^{(0)}, x^{(1)}, \dots , x^{(n)}\}$ in  $[a, b]_{\theta}$ such
that
\begin{equation}\label{eq9}
 a=x^{(0)}\stackrel{\theta}{<} x^{(1)}\stackrel{\theta}{<} \dots \stackrel{\theta}{<} x^{(n)}=b.
\end{equation}
If $P$ and $Q$ are two partitions of $[a, b]_{\theta}$ with $P\subseteq Q$, then $Q$ is called a 
{\bf refinement} of $P$.
\end{definition}

\medskip
Let  $[a, b]_{\theta}$ be a type $\theta$ closed interval. Suppose that 
$f: [a, b]_{\theta}\to \mathcal{R}^{(2)}$ is bounded and $P=\{x^{(0)}, x^{(1)}, \dots , x^{(n)}\}$ is a partion of  $[a, b]_{\theta}$. For $1\le i\le n$, the {\bf length} $\Delta x^{(i)}$ of the $i$-th type $\theta$ 
subinterval $[x^{(i-1)},\,x^{(i)}]_{\theta}$ is defined by $\Delta x^{(i)}:=x^{(i)}-x^{(i-1)}$. Clearly, 
 $\Delta x^{(i)}\stackrel{\theta}{>}0$, and  $\Delta x^{(i)}$ is  a real number if and only if 
$Ze(x^{(i)})=Ze(x^{(i-1)})$ for $1\le i\le n$. Since both $f_{Re}$ and $f_{Ze}$ are bounded on 
$[a, b]_{\theta}$, both $f_{Re}$ and $f_{Ze}$ are bounded on $[x^{(i-1)},\,x^{(i)}]_{\theta}$  for $1\le i\le n$. Hence, 
both 
$ \sup_if_{\clubsuit}:=\sup \big\{f_{\clubsuit}(x)\,|\, x\in [x^{(i-1)},\,x^{(i)}]\big\}$
and
$ \inf_if_{\clubsuit}:=\inf \big\{f_{\clubsuit}(x)\,|\, x\in [x^{(i-1)},\,x^{(i)}]\big\}$
exist as real numbers for $\clubsuit\in\{Re, \, Ze\}$ and  $1\le i\le n$. Based on these facts, 
we define the {\bf type $\theta$ upper sum} $U_{\theta}(P,  f)$ of $f$ with respect to the partition $P$ to be
\begin{eqnarray*}
U_{\theta}(P,  f)&=&\left\{\begin{array}{c}
\displaystyle\sum_{i=1}^n\big(\sup_if_{Re}+1^{\#}\,\sup_if_{Ze}\big)\Delta x^{(i)}\quad\mbox{$\theta =1$;}\\  \\
\displaystyle\sum_{i=1}^n\big(\sup_if_{Re}+1^{\#}\,\inf_if_{Ze}\big)\Delta x^{(i)}\quad\mbox{$\theta =2$}
\end{array}\right.
\end{eqnarray*}
and the {\bf type $\theta$ lower sum} $L_{\theta}(P,  f)$ of $f$ with respect to the partition $P$ to be
\begin{eqnarray*}
L_{\theta}(P,  f)&=&\left\{\begin{array}{c}
\displaystyle\sum_{i=1}^n\big(\inf_if_{Re}+1^{\#}\,\inf_if_{Ze}\big)\Delta x^{(i)}\quad\mbox{$\theta =1$;}\\  \\
\displaystyle\sum_{i=1}^n\big(\inf_if_{Re}+1^{\#}\,\sup_if_{Ze}\big)\Delta x^{(i)}\quad\mbox{$\theta =2$}.
\end{array}\right.
\end{eqnarray*}
Then the following four sets
\begin{equation}\label{eq10}
\big\{\clubsuit U_{\theta}(P,  f)\,\big|\, P\in \mathcal{P}_{\theta}\big\},\quad
\big\{\clubsuit L_{\theta}(P,  f)\,\big|\, P\in \mathcal{P}_{\theta}\big\}\quad\mbox{with}\quad
\clubsuit\in\{Re, \, Ze\}
 \end{equation}\label
are bounded subsets of the real number field $\mathcal{R}$, where $\mathcal{P}_{\theta}$ is the set of all partitions of  $[a, b]_{\theta}$, i.e.
$\mathcal{P}_{\theta}:=\big\{\, P\,\big|\, \mbox{$P$ is a partition of $[a,\,b]_{\theta}$}\,\big\}$. Hence, the supremums and infimums of the four sets in (\ref{eq10}) exist. Using these facts, we introduce the {\bf type $\theta$ lower integral} 
$\underline{\displaystyle\int_a^b} f(x)d_{\theta}x$ and  the {\bf type $\theta$ upper integral} 
$\overline{\displaystyle\int_a^b} f(x)d_{\theta}x$ of $f(x)$ on $[a, b]_{\theta}$  in the following way:
$$\underline{\displaystyle\int_a^b} f(x)d_{1}x=
sup\big\{Re L_{1}(P,  f)\,\big|\, P\in \mathcal{P}_{1}\big\}+
1^{\#}\,sup\big\{Ze L_{1}(P,  f)\,\big|\, P\in \mathcal{P}_{1}\big\},
$$
$$
\underline{\displaystyle\int_a^b} f(x)d_{2}x=
sup\big\{Re L_{2}(P,  f)\,\big|\, P\in \mathcal{P}_{2}\big\}+
1^{\#}\,inf\big\{Ze L_{2}(P,  f)\,\big|\, P\in \mathcal{P}_{2}\big\},
$$
$$
\overline{\displaystyle\int_a^b} f(x)d_{1}x=
inf\big\{Re U_{1}(P,  f)\,\big|\, P\in \mathcal{P}_{1}\big\}+
1^{\#}\,inf\big\{Ze U_{1}(P,  f)\,\big|\, P\in \mathcal{P}_{1}\big\}
$$
and
$$
\overline{\displaystyle\int_a^b} f(x)d_{2}x=
inf\big\{Re U_{2}(P,  f)\,\big|\, P\in \mathcal{P}_{2}\big\}+
1^{\#}\,sup\big\{Ze U_{2}(P,  f)\,\big|\, P\in \mathcal{P}_{2}\big\}.
$$
If the  type $\theta$ lower integral and the type $\theta$ upper integral of $f(x)$ on $[a, b]_{\theta}$ are equal, i.e.,
if $\underline{\displaystyle\int_a^b} f(x)d_{\theta}x=\overline{\displaystyle\int_a^b} f(x)d_{\theta}x$, then
we say that $f$ is {\bf type $\theta$ integrable} on $[a, b]_{\theta}$, we denote their common value by
$\displaystyle\int_a^b f(x)d_{\theta}x$ which is called the {\bf type $\theta$ integral} of $f$  on 
$[a, b]_{\theta}$.

\medskip
\begin{proposition}\label{pr3.1} Let $f: [a, b]_{\theta}\to \mathcal{R}^{(2)}$ is a bounded function, where $\theta\in\{1,\,2\}$.
\begin{description} 
\item[(i)] If $P$  and $P^*$ are partitions of $[a, b]_{\theta}$ and $P^*$ is a refinement of $P$, then
$$
 L_{\theta}(P,  f)\stackrel{\theta}{\le} L_{\theta}(P^*,  f)\stackrel{\theta}{\le} U_{\theta}(P^*,  f)
\stackrel{\theta}{\le} U_{\theta}(P,  f).
$$
\item[(ii)] $
\underline{\displaystyle\int_a^b} f(x)d_{\theta}x \stackrel{\theta}{\le}
\overline{\displaystyle\int_a^b} f(x)d_{\theta}x .
$
\item[(iii)] $f$ is  type $\theta$ integrable iff for each $\varepsilon\in \mathcal{R}^{(2)}$ with 
$\varepsilon\stackrel{\theta}{>}0$ and $(Re\,\varepsilon)(Ze\,\varepsilon)\ne 0$ there exists a  partition $P$ of  
$[a, b]_{\theta}$ such that
$$ U_{\theta}(P,  f)-L_{\theta}(P,  f)\stackrel{\theta}{<}\varepsilon .
$$\end{description}
\end{proposition}

\medskip
\noindent
{\bf Proof} The proof of Proposition \ref{pr3.1} follows from the definitions above and  the properties of the supremums and infimums. For example, let us prove (i) for $\theta=1$, i.e.,
\begin{equation}\label{eq11}
 L_{1}(P,  f)\stackrel{1}{\le} L_{1}(P^*,  f)\stackrel{1}{\le} U_{1}(P^*,  f)
\stackrel{1}{\le} U_{1}(P,  f).
\end{equation}

\medskip
The middle inequality  in (\ref{eq11}) follows directly from the definitions of type $1$ upper 
and lower sums. Suppose that 
$P=\{x^{(0)}, x^{(1)}, \dots , x^{(n)}\}$ and  consider the partition  $P^*$ formed by joining just  one 
point $x^*$ to $P$, where $ x^{(k-1)}\stackrel{1}{\le} x^*\stackrel{1}{\le} x^{(k)}$ for some
$k$ with $1\le k\le n$. Let
$$
\alpha_1 (f_{\clubsuit}):=inf \big\{\, f_{\clubsuit}(x)\,\big|\, x\in [x^{(k-1)}, x^*]_1\,\big\},
$$
$$
\alpha_2 (f_{\clubsuit}):=inf \big\{\, f_{\clubsuit}(x)\,\big|\, x\in [x^*, x^{(k)}]_1\,\big\},
$$
where $\clubsuit\in\{Re,\, Ze\}$. The terms in $L_{1}(P^*,  f)$ and $L_{1}(P,  f)$ are all the 
same except those over the subinterval $[x^{(k-1)},\,x^{(k)}]_1$. Thus we have
\begin{eqnarray}\label{eq12}
&&L_{1}(P^*,  f)-L_{1}(P,  f)\nonumber\\
&=&\big[(\alpha_1 (f_{Re})-\inf_kf_{Re})+1^{\#}(\alpha_1 (f_{Ze})-\inf_kf_{Ze})\big](x^*- x^{(k-1)})+
\nonumber\\
&&\quad +\big[(\alpha_2 (f_{Re})-\inf_kf_{Re})+1^{\#}(\alpha_2 (f_{Ze})-\inf_kf_{Ze})\big](x^{(k)}-x^*).
\end{eqnarray}
Since 
$$
\alpha_j (f_{\clubsuit})\ge \inf_kf_{\clubsuit}\quad\mbox{for $ j\in\{1,\, 2\}$ and  $ \clubsuit\in\{Re,\, Ze\}$ },
$$
we get
\begin{equation}\label{eq13}
(\alpha_j (f_{Re})-\inf_kf_{Re})+1^{\#}(\alpha_j (f_{Ze})-\inf_kf_{Ze})\stackrel{1}{\ge}0
\quad\mbox{for $ j\in\{1,\, 2\}$ }.
\end{equation}

Using  (\ref{eq13}) and the facts:
$x^*- x^{(k-1)}\stackrel{1}{>}0$ and $x^{(k)}-x^*\stackrel{1}{>}0$,
we get from  (\ref{eq12}) that
$L_{1}(P^*,  f)-L_{1}(P,  f)\stackrel{1}{\ge}0$ or
$L_{1}(P^*,  f)\stackrel{1}{\ge}L_{1}(P,  f)$.

\medskip
Similarly,  let
$$
\beta_1 (f_{\clubsuit}):=sup \big\{\, f_{\clubsuit}(x)\,\big|\, x\in [x^{(k-1)}, x^*]_1\,\big\},
$$
$$
\beta_2 (f_{\clubsuit}):=sup \big\{\, f_{\clubsuit}(x)\,\big|\, x\in [x^*, x^{(k)}]_1\,\big\},
$$
where $\clubsuit\in\{Re,\, Ze\}$. Since the terms in $U_{1}(P^*,  f)$ and $U_{1}(P,  f)$ are all the 
same except those over the subinterval $[x^{(k-1)},\,x^{(k)}]_1$, we have
\begin{eqnarray}\label{eq14}
&&U_{1}(P^*,  f)-U_{1}(P,  f)\nonumber\\
&=&\big[(\sup_kf_{Re}-\beta_2 (f_{Re}))+1^{\#}(\sup_kf_{Ze}-\beta_2 (f_{Ze}))\big](x^{(k)}-x^*)+
\nonumber\\
&&+\big[(\sup_kf_{Re}-\beta_1 (f_{Re}))+1^{\#}(\sup_kf_{Ze}-\beta_1 (f_{Ze}))\big](x^*- x^{(k-1)}).
\end{eqnarray}

Since
$$
\sup_kf_{\clubsuit}\ge \beta_j (f_{\clubsuit})\quad\mbox{for $ j\in\{1,\, 2\}$ and  $ \clubsuit\in\{Re,\, Ze\}$ },
$$
we have
\begin{equation}\label{eq15}
(\sup_kf_{Re}-\beta_j (f_{Re}))+1^{\#}(\sup_kf_{Ze}-\beta_j (f_{Ze}))\stackrel{1}{\ge}0
\quad\mbox{for $ j\in\{1,\, 2\}$}.
\end{equation}

It follows from  (\ref{eq14}) and  (\ref{eq15}) that
$U_{1}(P,  f)-U_{1}(P^*,  f)\stackrel{1}{\ge}0$ or $U_{1}(P,  f)\stackrel{1}{\ge}U_{1}(P^*,  f)$.

\medskip
This proves that (\ref{eq11}) holds.

\hfill\raisebox{1mm}{\framebox[2mm]{}}

\medskip
As a corollary of Proposition \ref{pr3.1} (iii), we have that if $f=f_{Re}+1^{\#}f_{Ze}$ is a function on a type $\theta$ closed interval  $[a, b]_{\theta}$ such that real-valued functions $f_{Re}$ and $f_{Ze}$ are continuous on the rectangle
$[Re\,a, Re\,b]\times [Ze\,2, \,Ze\,b]\subseteq \mathcal{R}^2=\mathcal{R}\times\mathcal{R}$, then $f$ is type $\theta$ integrable on $[a, b]_{\theta}$, where $\theta =1$ or $2$.

\medskip
The algebraic properties of the ordinary integral are still true for the type $\theta$ integrals.

\medskip
\begin{proposition}\label{pr3.2} Let $\theta =1$ or $2$ and let $k\in \mathcal{R}^{(2)}$ be a dual real number.
\begin{description} 
\item[(i)] If $f$ and $g$ are  type $\theta$ integrable on $[a, b]_{\theta}$, then
$f+g$ and $kf$ are  type $\theta$ integrable on $[a, b]_{\theta}$  and
$$\displaystyle\int_a^b (f+g) d_{\theta}x=\displaystyle\int_a^b f d_{\theta}x+
\displaystyle\int_a^b g d_{\theta}x,\quad
\displaystyle\int_a^b kf d_{\theta}x=k\displaystyle\int_a^b f d_{\theta}x .$$
\item[(ii)] If $f$ is type $\theta$ integrable on both 
$[a, c]_{\theta}$ and $[c, b]_{\theta}$, where $a\stackrel{\theta}{<}c\stackrel{\theta}{<}b$,
then $f$ is type $\theta$ integrable on both $[a, b]_{\theta}$ and 
$$
\displaystyle\int_a^b f\,d_{\theta}x=\displaystyle\int_a^c f\,d_{\theta}x
+\displaystyle\int_c^b f\,d_{\theta}x.
$$
\item[(iii)] If $f$, $g: [a, b]_{\theta}\to \mathcal{R}^{(2)}$ are type $\theta$ integrable and $f(x)\stackrel{\theta}{\ge} g(x)$ for all $x\in [a, b]_{\theta}$, then
$\displaystyle\int_a^b f(x)\,d_{\theta}x\stackrel{\theta}{\ge}\
\displaystyle\int_a^b g(x)\,d_{\theta}x$.
\end{description}
\end{proposition}

\medskip
\noindent 
{\bf Proof} Both (i) and (ii) are proved by using Proposition \ref{pr3.1}, and (iii) is proved by using the definitions of type $\theta$ inegrals and the properties of the two generalized order relations on the dual real number algebra.

\hfill\raisebox{1mm}{\framebox[2mm]{}}

\medskip
We finish this paper with the following counterpart of the odinary fundamental theorem of calculus in the context of dual real numbers.

\medskip
\begin{proposition}\label{pr3.3} Let $a$, $b\in \mathcal{R}^{(2)}$ and $a\stackrel{\theta}{<}b$, where $\theta=1$ or $2$.
\begin{description}
\item[(i)] If $f: [a, b]_{\theta}\to \mathcal{R}^{(2)}$ is  a function such that the real-valued functions  $f_{Re}$ and $f_{Ze}$ are continuous on the rectangle
$[Re\,a, \, Re\,b]\times [Ze\,a, \, Ze\,b]\subseteq \mathcal{R}^2$, then the function 
$F(x): [a, b]_{\theta}\to \mathcal{R}^{(2)}$ defined by
$$
F(x):=\displaystyle\int_a^x f(t)\,d_{\theta}t\quad\mbox{for $x\in [a, b]_{\theta}$}
$$
is type $\theta$ diffrential at  each $c\in [a, b]_{\theta}$ and $F'_{\theta}(c)=f(c)$.
\item[(ii)] If $f(x): [a, b]_{\theta}\to  \mathcal{R}^{(2)}$ is differential on $[a, b]_{\theta}$ and the derivative $f'(x)$ of $f(x)$ is integrable on $[a, b]_{\theta}$, then
$$\displaystyle\int_a^b f'(x)\,d_{\theta}x=f(b)-f(a).$$
\end{description}
\end{proposition}

\medskip
\noindent 
{\bf Proof} The way of proving Proposition \ref{pr3.3} comes from the application of the 
algebraic properties of type $\theta$ integrals in Proposition \ref{pr3.2}. Let us prove (i) to explain the way of doing the proofs.

\medskip
By the definitions of type $\theta$ integrals, we have
\begin{equation}\label{eq16}
a\stackrel{\theta}{\le} b\Longrightarrow
\displaystyle\int_a^b \,d_{\theta}x=b-a\quad\mbox{for $\theta\in\{1,\, 2\}$}.
\end{equation}

\medskip
Clearly, (i) holds if we can prove that for each positive real number $\epsilon>0$ there exists a  positive real number $\delta>0$ such that
\begin{equation}\label{eq17}
x\in N_{\theta}^*(c;\,\delta)\cap [a, b]_{\theta}\Rightarrow
\frac{||F(x)-F(c)-f(c)(x-c)||}{||x-c||}<\epsilon ,  
\end{equation}
where $\theta\in\{1,\, 2\}$. The proofs of (\ref{eq17}) for $\theta=1$ and $\theta=2$ are similar, so we prove (\ref{eq17}) for $\theta=1$. First,   
we choose two positive real numbers $\varepsilon_{Re}$ and $\varepsilon_{Ze}$ such that
\begin{equation}\label{eq18}
 \mbox{$0<\varepsilon_{Re}<\displaystyle\frac{\epsilon}{3}$\quad and\quad 
$0<\varepsilon_{Ze}<\displaystyle\frac{\epsilon -3\varepsilon_{Re}}{\sqrt{2}}$.}
\end{equation}

\medskip
Next, since both real-valued functions $f_{Re}$ and $f_{Ze}$ are continuous on the rectangle
$[Re\,a, \, Re\,b]\times [Ze\,a, \, Ze\,b]\subseteq \mathcal{R}^2$, both $f_{Re}$ and $f_{Ze}$ are uniformly 
continuous on the rectangle $[Re\,a, \, Re\,b]\times [Ze\,a, \, Ze\,b]$. Hence, there exist a positive real number $\delta >0$ such that
\begin{equation}\label{eq19}
\sqrt{(t_1-s_1)^2+(t_2-s_2)^2}<\delta\nonumber\\
\Rightarrow |f_{\clubsuit}(t_1,\,t_2)-f_{\clubsuit}(s_1,\,s_2)|<\varepsilon_{\clubsuit}
\end{equation}
for all $(t_1,\,t_2)$, $(s_1,\,s_2)\in [a_1, b_1]\times [a_2, b_2]$ and 
$\clubsuit\in\{Re,\, Ze\}$.

\medskip
Let $x\in N_{1}^*(c;\,\delta)\cap [a, b]_{\theta}$. Then $x\stackrel{1}{\ge}c$ and 
$x\stackrel{1}{\le}c$.

\medskip
{\bf Case 1}: $x\stackrel{1}{\ge}c$, in which case, by 
Proposition \ref{pr3.2} and (\ref{eq16}), we have
\begin{eqnarray}\label{eq20}
&&F(x)-F(c)-f(c)(x-c)=
\displaystyle\int_a^x f(t)\,d_{1}t-\displaystyle\int_a^c f(t)\,d_{1}t-f(c)(x-c)\nonumber\\
&=&\displaystyle\int_c^x f(t)\,d_{1}t-f(c)\displaystyle\int_c^x \,d_{1}t
=\displaystyle\int_c^x f(t)\,d_{1}t+\displaystyle\int_c^x \big(-f(c)\big)\,d_{1}t\nonumber\\
&=&\displaystyle\int_c^x \big[f(t)-f(c)\big]\,d_{1}t .
\end{eqnarray}

\medskip
For $c=c_1+c_21^{\#}\stackrel{1}{\le} t=t_1+t_21^{\#}\stackrel{1}{\le} x=x_1+x_21^{\#}$, where 
$c_1$, $c_2$, $t_1$, $t_2$, $x_1$ and  $x_2\in\mathcal{R}$, we have
\begin{eqnarray*}
&&\sqrt{(t_1-c_1)^2+(t_2-c_2)^2}\le \sqrt{(x_1-c_1)^2+(x_2-c_2)^2}\\
&\le& \sqrt{2(x_1-c_1)^2+(x_2-c_2)^2}=||x-c||<\delta\\
&\Rightarrow& |f_{\clubsuit}(t_1,\,t_2)-f_{\clubsuit}(c_1,\,c_2)|<\varepsilon_{\clubsuit}
\quad\mbox{for $\clubsuit\in\{Re,\, Ze\}$,}
\end{eqnarray*}
which implies
\begin{equation}\label{eq21}
-\varepsilon=-\varepsilon_{Re}-\varepsilon_{Ze}1^{\#}\stackrel{1}{\le}f(t)-f(c)
\stackrel{1}{\le}\varepsilon_{Re}+\varepsilon_{Ze}1^{\#}=\varepsilon.
\end{equation}

It folows from (\ref{eq16}), (\ref{eq21}) and Proposition \ref{pr3.2} (iii) that
$$
-\varepsilon (x-c)=\displaystyle\int_c^x (-\varepsilon)\,d_{1}t\stackrel{1}{\le}
\displaystyle\int_c^x \big[f(t)-f(c)\big]\,d_{1}t\stackrel{1}{\le}
\displaystyle\int_c^x \varepsilon\,d_{1}t=\varepsilon (x-c),
$$
which gives
\begin{equation}\label{eq22}
\Big|\clubsuit\Big(\displaystyle\int_c^x \big[f(t)-f(c)\big]\,d_{1}t\Big)\Big|
\le \big|\clubsuit\big(\varepsilon (x-c)\big)\big|\quad\mbox{for $\clubsuit\in\{Re,\, Ze\}$}.
\end{equation}
Since
\begin{eqnarray*}
&&\varepsilon (x-c)=(\varepsilon_{Re}+\varepsilon_{Ze}1^{\#})[(x_1-c_1)+(x_2-c_2)1^{\#}]\\
&=&\varepsilon_{Re}(x_1-c_1)+[\varepsilon_{Re}(x_2-c_2)+\varepsilon_{Ze}(x_1-c_1)]1^{\#}),
\end{eqnarray*}
we have
\begin{equation}\label{eq23}
Re\big(\varepsilon (x-c)\big)=\varepsilon_{Re}(x_1-c_1)
\end{equation}
and
\begin{equation}\label{eq24}
Ze\big(\varepsilon (x-c)\big)=\varepsilon_{Re}(x_2-c_2)+\varepsilon_{Ze}(x_1-c_1).
\end{equation}

\medskip
By (\ref{eq22}), (\ref{eq23}) and (\ref{eq24}), we have
\begin{eqnarray}\label{eq25}
&&\Big|\Big|\Big(\displaystyle\int_c^x \big[f(t)-f(c)\big]\,d_{1}t\Big)\Big|\Big|\nonumber\\
&=&\sqrt{2\Big[Re\Big(\displaystyle\int_c^x \big[f(t)-f(c)\big]\,d_{1}t\Big)\Big]^2
+\Big[Ze\Big(\displaystyle\int_c^x \big[f(t)-f(c)\big]\,d_{1}t\Big)\Big]^2}\nonumber\\
&\le&\sqrt{2}\Big|Re\Big(\displaystyle\int_c^x \big[f(t)-f(c)\big]\,d_{1}t\Big)\Big|
+\Big|Ze\Big(\displaystyle\int_c^x \big[f(t)-f(c)\big]\,d_{1}t\Big)\Big|\nonumber\\
&\le&\sqrt{2}\,\,\Big|\varepsilon_{Re}(x_1-c_1)\Big|+
\Big|\varepsilon_{Re}(x_2-c_2)+\varepsilon_{Ze}(x_1-c_1)\Big|\nonumber\\
&\le&\sqrt{2}\,\,\varepsilon_{Re}\big|x_1-c_1\big|+\varepsilon_{Re}\big|x_2-c_2\big|
+\varepsilon_{Ze}\big|x_1-c_1\big|
\end{eqnarray}

\medskip
It follows from (\ref{eq20}) and (\ref{eq25}) that
\begin{eqnarray}\label{eq26}
&&\frac{||F(x)-F(c)-f(c)(x-c)||}{||x-c||}
=\frac{\Big|\Big|\Big(\displaystyle\int_c^x \big[f(t)-f(c)\big]\,d_{1}t\Big)\Big|\Big|}
{||x-c||}\nonumber\\
&\le&
\frac{\sqrt{2}\,\,\varepsilon_{Re}\big|x_1-c_1\big|+\varepsilon_{Re}\big|x_2-c_2\big|+\varepsilon_{Ze}\big|x_1-c_1\big|}
{\sqrt{2(x_1-c_1)^2+(x_2-c_2)^2}}
\end{eqnarray}

\medskip
{\bf Case 11}: $x_1=c_1$, in which case, by (\ref{eq18}) and (\ref{eq26}), we get
\begin{equation}\label{eq27}
\frac{||F(x)-F(c)-f(c)(x-c)||}{||x-c||}
\le \frac{\varepsilon_{Re}\big|x_2-c_2\big|}{\sqrt{(x_2-c_2)^2}}=\varepsilon_{Re}<
\displaystyle\frac{\epsilon}{3}< \epsilon.
\end{equation}

\medskip
{\bf Case 12}: $x_1\ne c_1$ and 
$\left|\displaystyle\frac{x_2-c_2}{x_1-c_1}\right|\le \displaystyle\frac{1}{\sqrt{2}}$, in which case, using (\ref{eq18}) and (\ref{eq26}) again, we get
\begin{eqnarray}\label{eq28}
&&\frac{||F(x)-F(c)-f(c)(x-c)||}{||x-c||}\nonumber\\
&\le& 
\frac{\sqrt{2}\,\,\varepsilon_{Re}\big|x_1-c_1\big|+\varepsilon_{Re}\big|x_2-c_2\big|+\varepsilon_{Ze}\big|x_1-c_1\big|}
{\sqrt{2(x_1-c_1)^2+(x_2-c_2)^2}}\nonumber\\
&=&\frac{\sqrt{2}\,\,\varepsilon_{Re}+
\varepsilon_{Re}\left|\displaystyle\frac{x_2-c_2}{x_1-c_1}\right|+\varepsilon_{Ze}}
{\sqrt{2+\left(\displaystyle\frac{x_2-c_2}{x_1-c_1}\right)^2}}
<\sqrt{2}\,\,\varepsilon_{Re}+
\varepsilon_{Re}\left|\displaystyle\frac{x_2-c_2}{x_1-c_1}\right|+\varepsilon_{Ze}\nonumber\\
&<&\sqrt{2}\,\,\varepsilon_{Re}+\varepsilon_{Re}\displaystyle\frac{1}{\sqrt{2}}+\varepsilon_{Ze}
=\displaystyle\frac{3}{\sqrt{2}}\varepsilon_{Re}+\varepsilon_{Ze}
<\displaystyle\frac{1}{\sqrt{2}}\epsilon<\epsilon .
\end{eqnarray}

\medskip
{\bf Case 13}: $x_1\ne c_1$ and 
$\left|\displaystyle\frac{x_2-c_2}{x_1-c_1}\right|\ge \displaystyle\frac{1}{\sqrt{2}}$. In this case, we have  $\left|\displaystyle\frac{x_1-c_1}{x_2-c_2}\right|\le \sqrt{2}$ and
\begin{eqnarray}\label{eq29}
&&\frac{||F(x)-F(c)-f(c)(x-c)||}{||x-c||}\nonumber\\
&\le& 
\frac{\sqrt{2}\,\,\varepsilon_{Re}\big|x_1-c_1\big|+\varepsilon_{Re}\big|x_2-c_2\big|+\varepsilon_{Ze}\big|x_1-c_1\big|}
{\sqrt{2(x_1-c_1)^2+(x_2-c_2)^2}}\nonumber\\
&=&\frac{\sqrt{2}\,\,\varepsilon_{Re}\left|\displaystyle\frac{x_1-c_1}{x_2-c_2}\right|+
\varepsilon_{Re}+\varepsilon_{Ze}\left|\displaystyle\frac{x_1-c_1}{x_2-c_2}\right|}
{\sqrt{2\left(\displaystyle\frac{x_1-c_1}{x_2-c_2}\right)^2+1}}\nonumber\\
&\le& \sqrt{2}\,\,\varepsilon_{Re}\left|\displaystyle\frac{x_1-c_1}{x_2-c_2}\right|+
\varepsilon_{Re}+\varepsilon_{Ze}\left|\displaystyle\frac{x_1-c_1}{x_2-c_2}\right|
\nonumber\\
&\le& \sqrt{2}\,\,\varepsilon_{Re}\sqrt{2}+
\varepsilon_{Re}+\varepsilon_{Ze}\sqrt{2}
=3\varepsilon_{Re}+\varepsilon_{Ze}\sqrt{2}<\epsilon .
\end{eqnarray}

\medskip
It follows from (\ref{eq27}), (\ref{eq28}) and (\ref{eq29}) that (\ref{eq17}) holds in Case 1.

\bigskip
{\bf Case 2}: $x\stackrel{1}{\le}c$, in which case, a similar computation shows that (\ref{eq17}) holds in Case 2.

\bigskip
This proves that (\ref{eq17}) holds for each $x\in N_{1}^*(c;\,\delta)\cap [a, b]_{1}$.

\hfill\raisebox{1mm}{\framebox[2mm]{}}

\end{document}